 \documentclass[11pt]{article}
 \usepackage[top=1in,bottom=1in,left=1.5in,right=1.5in]{geometry}
  \usepackage{amsmath,amssymb}
  \usepackage{latexsym}
 \usepackage[dvips]{pstricks} 
  \usepackage{pst-node}
  \usepackage[dvips]{graphicx}

  \newcommand{\N}{\mathbb{N}}
  
  \newcommand{\R}{\mathbb{R}}

  \newcommand{\0}{\mathbf{0}}
  \newcommand{\1}{\mathbf{1}}

  \newcommand{\cA}{\mathcal{A}}

  \newcommand{\cM}{\mathcal{M}}

  \def\diag{\mathop{{\rm diag}}\nolimits}
  \newcommand{\hs}{\hspace*{\parindent}}
  \newcommand{\proof}{\hs \textbf{Proof.\ }}

  \newcommand{\trans}{^\top}
  \newcommand{\qed}{\hspace*{\fill} $\Box$\\}
  
  \newcommand{\haff}{\mathrm{haf\;}}
  \newcommand{\pfaf}{\mathrm{pfaf\;}}
  
  \newcommand{\perfmat}{\mathrm{perfmat\;}}

  \newtheorem{theo}{\bfseries \hs Theorem}[section]
  
  \newtheorem{prop}[theo]{\bfseries \hs Proposition}
  \newtheorem{lemma}[theo]{\bfseries \hs Lemma}
  
  \newtheorem{corol}[theo]{\bfseries \hs Corollary}

  \numberwithin{equation}{section} 

\begin{document}

 \title{Upper bounds for perfect matchings\\in pfaffian and planar graphs}

 \author
 {Afshin Behmaram\footnote{University of Tehran, Tehran, Iran, \emph{email}: behmaram@ut.ac.ir. This work was supported by ministry of science
of Iran and by UIC.}\; and
 Shmuel Friedland\footnote{Department of Mathematics, Statistics and Computer Science,
 University of Illinois at Chicago,
 Chicago, Illinois 60607-7045, USA,
 \emph{email}:friedlan@uic.edu. This work was supported by NSF grant DMS-1216393.}}

 \date{December 30, 2012}
 \maketitle

 \begin{abstract}
 We give upper bounds on weighted perfect matchings in pfaffian graphs.  These upper bounds are better than generalized Bregman's inequality.
 We show that some of our upper bounds are sharp for $3$ and $4$-regular pfaffian graphs.
 We apply our results to fullerene graphs.

 \end{abstract}



 \section{Introduction}

 The aim of this paper is give upper bounds on the number of matchings in pfaffian graphs using the Hadamard-Fischer determinant inequality.
 Let $G=(V,E)$ be a simple undirected graphs with the sets of $V$ vertices and $E$ edges.  Denote by $d(v)$ the degree of $v\in V$;
 $\Delta(G):=\max\{d(v), v\in V\}$; $\perfmat G$ the number of perfect matches in $G$; $K_n$ a complete graph on $n$ vertices; $C_n$ a cycle on $n$ vertices.
 The simplest result of this paper can be stated as follows.  Assume that $G$ is a planar graph, or more generally pfaffian.  Then
 \begin{equation}\label{upmbhadin}
 \perfmat G \le \prod_{v\in V} d(v)^{\frac{1}{4}}.
 \end{equation}
 We show that if $\Delta(G)\ge 3$ then the above inequality is stronger than the upper bound given in \cite{AF08}.
 Next we show that the above upper bound is sharp for the following $d$-regular connected planar graphs.  For $d=3$ two planar graphs: $K_4$ and $C_4\times K_2$, (the direct product of $C_4$ and $K_2$).  For $d=4$: a unique planar graph called octahedron, (a complement of a perfect match in $K_6$ and hence planar).   The rest of the paper is devoted to improvements of \eqref{upmbhadin}.

 We survey briefly the contents of this paper.
 In \S2 we recall the basic notions and results we use in this paper.
 In \S3 we apply Hadamard's determinant inequality to give an estimate on the sum of weighted
 perfect matches in a pfaffian graph with weighted edges.  We also give an upper bounds on perfect matchings in planar graphs with
 a planar girth $g'$, (which is not less than the girth of the graph.)  In \S4 we discuss examples where the inequality \eqref{upmbhadin}
 is sharp for $3$ and $4$-regular pfaffian graphs.  In \S5 we give upper bounds on the number of perfect matchings in a pfaffian graph $G=(V,E)$ with no $4$-cycles, in terms of the maximal match in the derived graph $G'=(V,E')$, where $(u,v)\in E'$ if and only if there is a path of length $2$ in $G$.
 In \S6.1, we introduce circular and semi circular graph and describe some classes of fullerene, (cubic graphs graphs with pentagon and hexagon faces), which have circular graph structure.  In \S 6.2 ,  the most technical section of this paper, we derive upper bounds of the number of perfect matches in semi-circular graphs with no $3$ and $4$-cycles.  As a corollary we show that a fullerene graph $G$ nontrivial cyclic 5-edge cut satisfies $\perfmat G\le 20^{\frac{|V|}{12}}$.

 \section{Basic notions and results}

 In this section we state our notation and recall some known results that will be used in this paper.
 Denote $n:=|V|$ the number of vertices; $m:=|E|$ the number of edges;
 $A(G)$ the adjacency matrix of $G$; $K_{r,r}$ a complete bipartite $r$-regular graph on
 $2r$ vertices.

 Let $G'=(V',E')$ be another graph.  Then the cartesian product of the graphs $G$ and $G'$, denoted by $G\times G'=(V\times V', \tilde E)$,
 is defined as follows. $((u,u'),(v,v'))\in \tilde E$ if and only if the following conditions hold.  Either $u=v$ and $(u',v')\in E'$ or
 $u'=v'$ and $(u,v)\in E$.

 Let $G=(V,E)$ be a simple graph.  Identify $V$ with $[n]:=\{1,\ldots,n\}$.
 Orient each edge $e\in E$.
 Then the skew symmetric adjacency matrix of
 $S(G)=[s_{ij}]_{i,j=1}^n$ is defined as follows. $s_{ij}=0$ if the undirected edge $(i,j)$ is not in $E$.  Assume that the undirected edge $(i,j)$ is
 in $E$.  Then $s_{ij}=1$ if the oriented edge goes from the vertex $i$ to vertex $j$.  Otherwise $s_{ij}=-1$.  Note that $A(G)=[|s_{ij}\emph{}|]$.

 Let $n$ be even an denote $n':=\frac{n}{2}$.
 A perfect match  $M$ in $G$ is a collection of $n'$ edges in $G$  which do not have a common vertex.  Denote by $\cM(G)$ the set of perfect matches
 in $G$.  Then $\perfmat G=|\cM(G)|$.

 Let $w: E\to \R_+:=[0,\infty)$ be the weights of the edges $E$.  Denote $G_w:=(V,E,w)$ the weighted
 graph.  Then $A(G_w)=[w(u,v)]_{u,v\in V}$ is the weighted adjacency matrix of $G_w$, where $w(u,v)=0$ if the edge $(u,v)\not\in E$.
 Denote by $d_w(v):=\sum_{u,(u,v)\in E} w(u,v)$ the weighted degree of $v\in V$.
 Direct the edges of $G$.  Then the corresponding weighted sign adjacency matrix is $S(G_w)=[s(u,v)]_{u,v\in V}$
 is a skew symmetric matrix, where $s(u,v)\ge 0$ if the oriented edge is from $u$ to $v$.

 Let $M\in\cM(K_n)$.  Then $M=\{(i_1,j_1),\ldots,(i_{n'},j_{n'})\}$ where $\{i_1,j_1,\ldots,i_{n'},j_{n'}\}=[n]$.
 We will assume that $i_k<j_k$ for $k\in[n']$ and $1\le i_1<i_2<\ldots<i_{n'}\le n$.  Associate with $M$ the following permutation
 $\sigma_M:[n]\to [n]$: $1\to i_1,2\to j_1,\ldots,n-1\to i_{n'},n \to j_{n'}$.

 Let $W=[w_{ij}]_{i,j=1}^n$ be a symmetric matrix with zero diagonal.   The \emph{haffnian} of $W$ is defined as
 \begin{equation}\label{defhaf}
 \haff W = \sum_{\{(i_1,j_1),\ldots,(i_{n'}j_{n'})\}\in \cM(K_n)} w_{i_1j_1}\ldots w_{i_{n'}j_{n'}}.
 \end{equation}

 Clearly, the entries of $W$ can be identified with the weights of edges in $K_n$.  Then $\perfmat (K_n)_w:=\haff W$ is the sum of
 all weighted perfect matches of $(K_n)_w$.  Assume that $w_{ij}=0$ if $(i,j)$ is not an edge in $E$.  Then $W$ represents $G_w$
 and $\perfmat G_w:=\haff W$ is the sum of all weighted perfect matches of $G_w$.

 Let $S=[s_{ij}]_{i,j=1}^n$ be a skew symmetric matrix.  Recall the definition of the pfaffian of $S$:
 \begin{equation}\label{defpfaf}
 \pfaf S = \sum_{\{(i_1,j_1),\ldots,(i_{n'}j_{n'})\}\in \cM(K_n)}\textrm{sign}(\sigma_M) s_{i_1j_1}\ldots s_{i_{n'}j_{n'}}.
 \end{equation}
 Recall the well known formula $\det S=(\pfaf S)^2$.  We now state a well known result and bring its proof for completeness.
 \begin{lemma}\label{hafdetin}  Let $G=(V,E)$ be a simple graph on even number of vertices.  Let $w: E\to (0,\infty)$.
 Then
 \begin{equation}\label{detpmbd}
 \det S(G_w)\le (\perfmat G_w)^2.
 \end{equation}
 Equality holds if and only if $\textrm{sign}(\sigma_M)= \textrm{sign}(\sigma_{M'})$  for any  $M,M'\in \cM(G)$.
 \end{lemma}
 \proof  Clearly $|\pfaf S(G_w)|\le \perfmat G_w$.  Since each edge in $G$ has a positive weight it follows that equality holds if and
 only if $\textrm{sign}(\sigma_M)= \textrm{sign}(\sigma_{M'})$  for any  $M,M'\in \cM(G)$.
 As $\det S(G_w)=(\pfaf S(G_w))^2$ we deduce the lemma.  \qed

 A graph is called \emph{pfaffian} if there is an orientation of $E$ such that $\det S(G)=(\perfmat G)^2$.
 We call such orientation a pfaffian orientation.
 It was shown by Kasteleyn \cite{Kas67} that every planar graph is pfaffian.
 Recall that an orientation of $E$ is pfaffian if and only if the following condition hold.
 Let $C$ be an even cycle of $G$.  If the subgraph $G\setminus C$ has a perfect matching then
 $C$ has an odd number of edges directed in either direction of the cycle.
 An example of a pfaffian nonplanar graph is $K_4\times K_2$.  See \cite{Tho06} for a recent review on pfaffian graphs.
 Assume that $G$ is pfaffian.  It is straightforward to see that for a pfaffian orientation and any $w:E\to [0,\infty)$ we have equality in \eqref{detpmbd}.

 Let $B=[b_{uv}]_{u,v\in V}$ be a real symmetric nonnegative definite matrix of order $|V|$, which is denoted by $B\succeq 0$.
 The generalized Hadamard-Fischer inequality, abbreviated here as H-F inequality, states
 \begin{equation}\label{FHineq}
 \det B[U\cup W]\det B[U\cap W]\le \det B[U]\det B[W] \textrm{ where } U,W\subset V.
 \end{equation}
 Here $\det B[\emptyset]=1$.  (See for example \cite{FG12} and references therein.)
 Assume that $U\cap W=\emptyset$.
 Furthermore if the left hand-side of \eqref{FHineq} is positive then equality holds if and only if $B[U\cup W]$ is a block diagonal matrix $\diag(B[U],B[W])$.
 The above inequality for $U\cap W=\emptyset$ is equivalent to the Hadamard-Fischer inequality.
 \begin{equation}\label{hadfisin}
 \det B\le \prod_{i=1}^l \det B[U_i],
 \end{equation}
 where $\cup_{i=1}^l U_l$ is a partition $V$ to a union of nonempty disjoint sets.
 In particular Hadamard's determinat inequality for $B\succeq 0$ is
 \begin{equation}\label{hadineq}
 \det B\le \prod_{v\in V} b_{vv}.
 \end{equation}
 If all diagonal entries of $B$ are positive then equality holds if and only if $B$ is a diagonal matrix.

 We apply our result for fullerene graphs.  Fullerenes are famous chemical graphs that were introduced by Harold Kroto, Richard Smalley and Robert Curl  in 1985 and in 1996.  The Nobel prize of Chemistry was awarded to them for discovery of fullerenes.  Fullerene graph is a planar cubic graph whose faces are pentagons or hexagons.  Let $ F_n$ be a fullerne with n vertices.  The smallest fullerene is $ F_{20}$ is called dodecahedron, which consist of twelve pentagons and no hexagon.  Let p and h be the number of pentagons and hexagons in fullerene  graph with n vertices and m edges. By Euler's formula we have :  $ p=12,   n=2h+20,   m=3h+30$.

 \section{Upper bounds on perfect matchings in pfaffian graphs}

 \begin{theo}\label{uppefmatbd} Assume that $n$ is even, $G$ has a pfaffian orientation and $w:E\to [0,\infty)$.
 Then
 \begin{equation}\label{uppefmatbd1}
 \perfmat G_w\le \prod_{v\in V} (d_{w,2}(v))^{\frac{1}{4}}, \quad d_{w,2}(v):=\sum_{u, (u,v)\in E} w(u,v)^2.
 \end{equation}
 In particular \eqref{upmbhadin} holds.
 \end{theo}
 \proof
 Assume that $B=S(G_w)S(G_w)\trans=-S(G_w)^2$.
 Note that $B[\{v\}]=d_{w,2}(v)$.  Hadamard's determinant inequality yields that   $\det B \le \prod_{v\in V} d_{w,2}(v)$  and $\det B=\det(S(G))^2=(\perfmat G_w)^4$  . This relation yields the theorem. \qed

 We next recall the upper bound on then number of perfect matchings in $G$ given in \cite{AF08},
 which was mentioned in the Abstract as generalized Bregman's inequality.
 \begin{equation}\label{AFin}
 \perfmat G\le \prod_{v\in V} (d(v)!)^{\frac{1}{2d(v)}}.
 \end{equation}
 Equality holds if and only if $G$ is a union of complete regular bipartite graphs.

 Clearly $(d!)^{\frac{1}{2d}}=d^{\frac{1}{4}}$  for $d=1,2$ .  For an integer $d\ge 3$ we have the following inequality.
 \begin{lemma}\label{dgreat3in}
 \begin{equation}\label{basineq}
 (d!)^{\frac{1}{2d}} > d^{\frac{1}{4}} \quad \textrm{ for } d=3,\ldots.
 \end{equation}
 \end{lemma}
 \proof  Stirling's formula for $d!$ implies the inequality
 $(d!)^{\frac{1}{d}}\ge \frac{d}{e}$.  Next observe that  $\sqrt{d}< \frac{d}{e}$  for  $d\ge 8$.
 One easily verifies that \eqref{basineq} holds for $d=3,\ldots,7$.
 \qed

 As $\perfmat K_{r,r}=r!$ we deduce a well known result.
 \begin{corol}\label{krrnpf}  $K_{r,r}$ is not pfaffian for $r\ge 3$.
 \end{corol}

 Recall that the girth $g(G)$ of a graph $G$ is the length of the shortest cycle in $G$.
 Let $G$ be a planar connected graph which is not a tree.  We define a planar girth of $G$ as follows.
 Assume that we have a realization of $G$ as a partition of the plane to $f$ faces $F_1,\ldots,F_f$ including the unbounded face.
 Denote by $\ell(F_i)$ the length of the face $F_i$.  Then the maximum value of $\min \ell(F_i)$ for all possible realization of $G$ as a planar
 graph is called the \emph{planar girth} of $G$, and is denoted by $g'(G)$.
 A realization of $G$ where $g'(G)$ is the minimal length of the faces is called an optimal realization.
 It is easy to construct an example of a connected planar graph $G$
 with $g(G)=3$ and $g'(G)>3$.  Let $g'(T)=\infty$ for a tree. For a nonconnected planar graph $G$ we define $g'(G)$ to be the minimum of the planar
 girths of its connected components.

 The following result is probably well known and we bring its proof for completeness.
 \begin{lemma}\label{gpineq}
 Let $G$ be a planar graph with a finite planar girth $g'$.  Then
 \begin{equation}\label{mgnin}
 m \le \frac{g'}{g'-2}(n-2).
 \end{equation}
 Equality holds if and only if $G$ is connected, and there exists an optimal realization of $G$, such that each face, including the outside face, is bounded by a cycle of length $g'\ge 3$.
 \end{lemma}
 \proof
 Assume that $G$ is connected.  Let $F_1,\ldots,F_f$ be the $f$ faces of an optimal realization of $G$.  Euler's formula states $f-m+n=2$.  Clearly
 $2m=\sum_{i=1}^f g(F_i)\ge g'f=g'(2+m-n)$, which implies \eqref{mgnin}.
 Assume that $G$ has $l$ connected components $G_i,i=1,\ldots,l$.   Then $m_i \le \frac{g'}{g'-2}(n_i-2)$ if $G_i$ has a cycle.
 If $G_i$ has no cycle, i.e. $G_i$ is a tree, tbmhen $m_i=n_i-1$.  This implies \eqref{mgnin}.
 The equality case follows straightforward.  \qed

 \begin{theo}\label{upbdg'} In every planar graphs with a finite planar girth $g'$ the number of perfect matchings are less than
 $(\frac{2g'}{g'-2})^{\frac{n}{4}}$.  In particular, if $g'>3$, i.e. there is a realization of $G$ such that each face is not a triangle,  then
 the number of perfect matchings is less than $2^{\frac{n}{2}}$.
 \end{theo}

 \proof
 Apply the AM�-GM inequality and the fact that the sequence $(1-\frac{2}{l})^{\frac{l}{2}}$ is an increasing sequence
 converging to $e^{-\frac{1}{2}}$ to deduce
 \begin{equation}\label{uppefmatbd2}
  \perfmat G\le \prod_{v\in V} (d(v))^{\frac{1}{4}}\le(\frac{2m}{n})^\frac{n}{4}\le (\frac{2g'(n-2)}{(g'-2)n}) ^{\frac{n}{4}}<e^{-\frac{1}{2}}(\frac{2g'}{g'-2})^{\frac{n}{4}}.
 \end{equation}
 If $g'\ge 4$ then $\frac{2g'}{g'-2}\le 4$ and the above inequality implies that
 $\perfmat G< 2^{\frac{n}{2}}$.
 \qed

 \section{Sharp upper bounds on $3$ and $4$-regular pfaffian graphs}

 Let $\phi(k,G)$ be the number of matches which cover $2k$ vertices of $G$.   Define $G\times K_2$, the cartesian product of $G$ and $K_2$, as two copies of $G$ where we join by an edge two copies of each vertex $v\in V$.  It is easy to show that $C_4 \times K_2$ is a planar graph.
 Clearly $\perfmat C_4=2, \perfmat K_4=3$.  Hence
 \begin{eqnarray}\label{pmatC4K2}
 &&\perfmat(C_4\times K_2)=\perfmat(C_4)^2 +\phi(1,C_4)+1=4+4+1=3^2,\\
 &&\label{pmatK4K2}
 \perfmat(K_4\times K_2)=\perfmat(K_4)^2+\phi(1,K_4)+1=9+6+1=4^2.
 \end{eqnarray}
 It is well known that $K_4\times K_2$ is pfaffian \cite{Tho06}.
 \begin{theo}\label{3regupbounds}  Let $G=(V,E)$ be a $3$-regular simple pfaffian graph.  Then the number of perfect matchings
 in $G$ is not more than $3^{\frac{|V|}{4}}$.  For a planar $3$-regular graph $G$ equality holds if and only if $G$ is a disjoint union of copies of $K_4$ or $C_4\times K_2$.
 \end{theo}
 \proof  Hadamard's determinant inequality gives the upper bound.  To discuss the equality case for a planar graph $G$ it is enough to show that each connected component $H$ of $G$ is either $K_4$ or $C_4\times K_2$.  Since the graph is simple, each $H$ contains at least $4$ vertices.
 If $H$ has $4$ vertices then $H=K_4$.
 Assume to that $H=(V_1,E_1)$ contains more than $4$ vertices.  Since we assumed that $\det S(H) =3 ^{\frac{|V_1|}{2}}$
 it follows that $S(H)^2=-3I$, where $I$ is the identity matrix.  This means that if there is a path of length two connecting two vertices $u,v\in V_1$  then we have exactly another path from $u$ to $v$, as $H$ is $3$-regular.
 Consider the faces generated by the planar graph $H$.  Suppose first that we have a face bounded by a triangle with vertices $u,v,w$.
 Since there is a path of length two between any two vertices of this triangle,
 each edge of this triangle must be bounded by another triangle.  The assumption $H$ is $3$-regular implies that $H=K_4$, contrary that $H$ has more than $4$ vertices.

 So $g'(H)\ge 4$.  Suppose there is a face $F_1$ of $H$ of which is bounded by a cycle of length $5$ at least.  Take three consecutive points
 $u_1-u_2-u_3$ on this cycle.  Since the distance from $u_1$ to $u_3$ is $2$ $H$ contains another path of length two between $u_1$ and $u_2$,
 namely $u_1-v-u_3$.  This path must lie on one or two adjacent faces of $H$ to $F_1$.  Since the degree of $u_2$ is $3$
 this path lies on two adjacent triangles with vertices $u_1,u_2,v$ and $u_2,v,u_3$.  This contradicts our claim that $g'(H)>3$.
 Hence each face of $H$ is bounded by a $4$-cycle.  Take one compact face $F_1$ bounded by the $4$-cycle $u_1-u_2-u_3-u_4-u_1$.  Since the degree of each vertex of $H$ is $3$ there must be another adjacent face of $F_1$ bounded by the cycle $u_1-u_2-v_1-v_2-u_1$, call it $F_2$.  Since the distance between $v_1$ and $u_3$ is $2$ we must have a $4$-cycle $u_3-u_2-v_1-v_3-u_3$ which bounds the face $F_3$ adjacent to $F_1$ and $F_2$.
 Similarly, we must have a $4$-cycle $u_4-u_1-v_3-v_4$  which bounds the face $F_4$ adjacent to $F_1$ and $F_2$.  Note that $v_3\ne v_4$.
 Since the dist$(v_2,v_3)=$dist$(v_1,v_4)=2$ we must have an edge $v_3-v_4$ in $H$.  Then $H=C_4\times K_2$.  \qed

 The octahedron graph, denoted by $O$, is a $4$-regular graph with $6$ vertices.  It is obtained from $K_6$ by deleting $3$ edges forming a complete match.
 Hence $O$ is unique up to an isomorphism.
 It is easy to see that $O$ is planar.  A straightforward calculation shows that $\perfmat O=8$.
 Indeed, recall that $\perfmat K_6=15$.  Take a perfect match $M$ in $K_6$ and count all matches in $K_6$ which contain
 at least one edge from $M$.  One match is $M$.  All other matches contain exactly one edge from $M$.  Let $(u,v)\in M$.
 Then $3=\perfmat K_4$ is the number of matches in $K_6$ containing $(u,v)\in M$.  One of them is $M$.  Hence the number of all
 matches in $K_6$ which contain $(u,v)$ and no other edges from $M$ is $2$.  Thus $\perfmat O=15 -(1+3 \times 2)=8$.

 Let $G$ be a connected $4$-regular planar graph with more than $6$ vertices.
 In \cite[\S2]{Leh} it is shown that every 4-regular connected planar graph $G$ on more than $6$ vertices
 contains one of the configurations A, B, C, or D.  In particular, it follows that not all the faces of $G$ are triangles.

 \begin{theo}\label{4regupbounds}  Let $G=(V,E)$ be a $4$-regular simple pfaffian graph on an even number of vertices.  Then the number of perfect matchings is not more than $4^{\frac{|V|}{4}}$.  For planar graphs equality holds if and only if $G$ is a disjoint union of copies the octahedron graph.
 \end{theo}
 \proof
 The inequality $\perfmat G\le 2^{\frac{n}{2}}$ follows from Hadamard's determinant inequality.  Equality holds if and only if $S(G)^2=-4I$.
 As the number of perfect matches in $O$ is $2$ to the power of a half of the vertices we deduce that $S(O)^2=-4I_6$.

 Let $G$ be a $4$-regular planar graph which satisfies $S(G)^2=-4I$.
 The inequality \eqref{detpmbd} implies that the orientation inducing $S(G)$ is pfaffian.
 Hence each connected component $H$ of $G$ satisfies $S(H)^2=-4I$.
 Furthermore $S(H)$ is a pfaffian orientation of $H$ and $\perfmat H=2^{\frac{|V_1|}{2}}$.
 Assume to the contrary that $H=(V_1,E_1)$ is not $O$, i.e. $|V_1|\ge 8$.

  Observe next that the equality $S(H)^2=-4I$ implies  the following condition.  Assume that there is a path of length two connecting two vertices $u,v\in V_1$.   Then we have either one or three additional paths of length two from $u$ to $v$, as $H$ is $4$-regular.
 Theorem \ref{upbdg'} implies that $g'(H)=3$.  Take a face $F_1$ of $H$ which is a triangle which vertices  $v_1,v_2,v_3$.
 As in the proof of Theorem \ref{3regupbounds} each neighboring face  $F_1$ is a triangle.   Hence all faces of $H$ are triangles.
 This contradicts the theorem of \cite{Leh}.
 \qed
the example of nonplanar graph which equality hold is the disjoint uion of  $ K_4\times K_2 $.

 \section{Applications of H-F inequality to pfaffian graphs}
 \begin{theo}\label{no4cyc}  Let $G=(V,E)$ be a pfaffian connected graph.  Denote by $G'=(V,E')$ be the induced graph by $G$, where $(u,v)\in E'$
 if and only if there exists a path of length $2$ between $u$ and $v$ in $G$.  Let $M'$ be a match in $G'$.  Assume that $G$ does not have $4$-cycles.
 Then
 \begin{equation}\label{no4cyc1}
 \perfmat G\le (\prod_{(u,v)\in M'} (d(u)d(v)-1))^{\frac{1}{4}}(\prod_{v\in V\setminus{V(M')}}d(v))^{\frac{1}{4}}.
 \end{equation}
 \end{theo}
 \proof  Consider the matrix $B=-S(G)^2\succeq 0$.  Let  $B=-S(G)^2=[b_{u v}]_{u,v\in V}$.  Assume that $u\ne v$.
 The assumption that $B$ has no  $4$-cycles
 means that $b_{uv}=\pm 1$  if there exists a path of length $2$ between $u$ and $v$ and otherwise $b_{uv}=0$.  Clearly, $b_{uu}=d(u)$.
 Associate with $M'$ the following partition $[n]$.
 With each pair $(u,v)\in M'$ associate $J_{u,v}:=\{u,v\}$.  For all vertices $w\in V$ not covered by $M'$ we associate $J_w=\{w\}$.
 Apply \eqref{hadfisin} to deduce \eqref{no4cyc1}.  \qed

 One can use the following proposition to estimate $|M'|$.
 \begin{prop}\label{estM'}
 Assume that $G=(V,E), n=|V|$ is a connected graph with a path of length $l\ge 3$.  Then the induced  graph $G'$ given in Theorem \ref{no4cyc}
 has a match $M'$ of cardinality $2\lfloor \frac{l}{4}\rfloor$ and $\lfloor \frac{l}{2}\rfloor$ if $l$ is even or odd respectively.
 \end{prop}

 It is conjectured that any fullerene graph has a hamiltonian cycle.  It is known that any fullerene graph has a cycle of length
  $\frac{5n-4}{6}$ at least \cite{KPSS09}.
 \begin{corol}\label{estM'ful}
 Assume that $G=(V,E)$ is a fullerene graph.  Then $\perfmat G$ is bounded above by the following quantities.
 \begin{equation}\label{estM'a}
 \begin{cases} 8^{\frac{n}{8}}   \text{ if } G \text{ has a hamiltonian cycle and } n \text{ is divisible by } 4
 \\
 8^{\frac{n-2}{8}} \sqrt{3}   \text{ if } G \text{ has a hamiltonian cycle and } n \text{ is not divisible by } 4
 \\
 8^{\frac{1}{4}\lfloor \frac{5n-4}{12}\rfloor}3^{\frac{n-2\lfloor \frac{5n-4}{12}\rfloor}{4}}.
 \end{cases}
 \end{equation}

 \end{corol}

 \begin{theo}\label{upboundfullerene} Every cubic pfaffian graph with n vertices, which has a perfect match and no $4$-cycles, has at most $8^{\frac{n}{12}}3^{\frac{n}{12}}$ perfect matching.
 \end{theo}
 \proof Let $M$ be a perfect match in $G$.  Then the complement of $M$ in $G$ is a union of spanning cycles $C_{p_1},\ldots,C_{p_k}$.
 Each even cycle $C_{p_i}$ induces a match of order $2\lfloor\frac{p_i}{4}\rfloor$ in the induced graph $G'$.  Each odd cycle $C_{p_j}$ induces
 a match of order $\frac{p_j-1}{2} $ in the induced graph $G'$.  Let $M'$ be a match in $G'$ induced as above by the cycles $C_{p_1},\ldots,C_{p_k}$.  The smallest $M'$ will be obtained when each $p_i=3$.  So $|M'|\ge \frac{n}{3}$.  Apply Theorem
 \ref{no4cyc1} to deduce the theorem.  \qed

 \section{Circular and semi-circular graphs}\label{sub:circsemi}
 \subsection{Definition and properties}\label{subsub:defprop}
 A connected planar graph $G=(V,E)$ is called \emph{semi-circular} if $G$ can be obtained as follows.  Draw $k\ge 1$ circles $O_1,\ldots,O_k$ in the plane of radii $0<r_1<\ldots<r_k$ centered at the origin.  At each circle $j$ we assume that we have $m_j\ge 3$ vertices for $j=1,\ldots,k$, denoted as $U_j \subset V$.  The subgraph of $G(U_j)$ is the cycle $C_{m_j}$ given by $O_j$ and the vertices $U_j$ for $j\in [k]$.
 In addition connect some vertices in two adjacent circles $U_i$ and $U_{i+1}$ such that the corresponding edges lying in the annulus
 $r_i\le \sqrt{x^2+y^2}\le r_{i+1}$ and do not intersect in the interior $r_i< \sqrt{x^2+y^2}< r_{i+1}$, for $i=1,\ldots,k-1$.
 The rest of the vertices $V_0$ of $G$ lie inside the open disk $x^2+y^2< r_1$.   So $G(U_1\cup V_0)$ is a connected planar graph
 whose vertices and edges lie in the closed disk $x^2+y^2\le r_1$.
 (Since $G$ is connected we assume that we have at least one edge connecting the sets $V_i$ and $V_{i+1}$ for $i\in [k-1]$.)
 A semi-circular graph is called circular if  $V_0=\emptyset$.

 The following characterization of a semi-circular graph is straightforward.
 \begin{lemma}\label{charcircg}  A connected planar graph $G$ is semi-circular if and and only if it can be realized in the plane as follows.
 Assume that $G=(V,E)$ induces $f-1$ bounded faces $F_1,\ldots,F_{f-1}$ and an unbounded face $F_f$ in $\R^2$.  Then there exists a following sequence of subgraphs of $G$: $G=G_k=(V_k,E_k)\supsetneqq G_{k-1}=(V_{k-1},E_{k-1})\supsetneqq\ldots\supsetneqq G_1=(V_1,E_1)$ with the following properties.  The bounded faces of each $G_i$ are
 $ F_j, j\in\cA_i$, where $\cA_k=[f-1]\supsetneqq \cA_{k-1}\supsetneqq\ldots\supsetneqq\cA_1\ne \emptyset$.  The unbounded face of $G_i$ is a cycle
 $C_{m_i}$, where $m_i\ge 3$.   $V_{i-1}$ are all the vertices of $V$ contained in the interior of $\cup_{j\in\cA_i} F_j$, for $j=k,\ldots,1$.
 Furthermore $G$ is circular if $V_0=\emptyset$.

 \end{lemma}

 Clearly, any connected planar graph with edges and no bridges is a semi-circular graph with $k=1$.
 So a really interesting semi-circular graphs is a graph where $V_0$ is a relatively small set in $V$.
 Available representations of $26,70$ fullerene graphs are only semi-circular with a relative small $V_0$.
 The $20, 24, 60$ fullerene graphs are circular graphs.

 We now define two infinite families of circular fullerene graphs see  \cite[8-4]{LoeblL8}.
 The first circle and the last circle is $C_k$ with $k=5$ or $k=6$.  I.e. each circle bounds a bounded and unbounded either pentagon or hexagon face respectively.   All other $\ell$ circles, $\ell\ge 1$,  are $C_{2k}$.  That is, the first and the last $C_k$ are surrounded by layers of $k$ pentagons.
 Note that the number of vertices in such fullerene graph is $2k(\ell+1)$.
 We call each family of these fullerene graphs pentacap and hexacap for $k=5,6$ respectively.
 In  \cite{KS08,KM08} it was proved that a fullerene graph with a non trivial 5-edge cut is a pentacap with $l\ge 2$.

 We observe that we can slightly generalize the construction of pentacaps and hexacaps fullerenes.
 \begin{prop}
 Let $F=(V,E)$ be a circular fullerene graph where the first circle is $C_k$ and the second circle is $C_{2k}$ for $k=5,6$.
 (I.e. the first circle is surrounded by a layer of $k$ pentagons.)
 Then there exists a circular fullerene graph $F'=(V',E')$ with $|V'|=|V|+k+1$.
 \end{prop}
 \proof   Add $k$ vertices, $u_1,\ldots,u_k$ to the first circle $C_k$, one new between two old.  Add a new inner circle $C_k$ with vertices $w_1,\ldots,w_k$
 and connect $u_i$ with $w_i$ for $i=1,\ldots,k$.  \qed

 Recall the leapfrog operation on a simple connected bridgeless planar graph $G=(V,E)$ with faces $F_1,\ldots,F_f$, including the unbounded face.
 Choose a point $w_i$ in the interior of each $F_i$.  Let $W=\{w_1,\ldots,w_{f}\}$.  Then $\tilde G=(\tilde V,\tilde E)$ where $\tilde V=V\cup W$ and
 $\tilde E=E\cup E_1$.  The edges $E_1$ are obtained  by connecting each $w_i$ to each vertex $V$ on the boundary of $F_i$.
 Since $G$ is bridgeless it follows that $|E_1|=2|V|$.
 Note that all faces $\tilde F_1,\ldots, \tilde F_{\tilde f}$ of $\tilde G$ are triangles.  Then $Le(G)=(V',E')$ is a  connected planar graph, where $|V'|=|\tilde f|$.  Each vertex $v_i'\in V'$ is an interior point of $\tilde F_i$.  $(v_i',v_j')\in E'$ if and only if $\tilde F_i$ and $\tilde F_j'$ have a common edge.
 So $Le(G)$ is the Poincar\'e dual of $\tilde G$.  It is known that if $G$ is fullerene then $Le(G)$ is fullerene \cite{FR96}.

 \begin{theo}\label{leapfrog}
 Let $F$ be a pentacap or hexacap.  That is  $F$ is a circular fullerene with inner and outer circles $C_{k}, k=5,6$ and $\ell$ circles $C_{2k}$.  Then $Le(F)$ is a circular fullerene if and only if $\ell$ is odd.
 Furthermore, if $\ell$ is odd then the circular structure of $Le(F)$ is as follows.  The first inner circle is $C_k$, the second circle is $C_{3k}$, then there are $\frac{3\ell-1}{2}$ circles $C_{4k}$, then follows the circle $C_{3k}$, and the last circle is $C_k$.
 If $\ell$ is even then $Le(F)$ has the following semi-circular structure.  The first circle is $C_{k}$ then the next  $\frac{3l}{2}$ circles are $C_{4k}$.  The
 circle number $\frac{3\ell}{2}+1$ is $C_{3k}$.  The last circle is $C_{k}$.  Inside the first circle there are $k$ vertices.
 \end{theo}

 \proof  We identify the vertices of $Le(F)$ with the triangle faces of $\tilde F$.
 The $k$ triangles inside the first circle $C_k$ are denoted by $\alpha_1$ and the $k$ triangles outside the last circle of $C_k$ of $F$ are denoted by
 $\alpha_2$.  Consider first the $k$ pentagons surrounding the first circle in $F$.  Each such a pentagon is divided to
 $5$ triangles in $\tilde F$.  The triangle which has one edge on the first circle we name $a_1$.  The two neighboring triangles of $a_1$ in this pentagon
 are called $b_1$.  The other two triangles of this pentagon are called $c_1$.  Each $c_1$ has one edge on the second circle of $F$.
 The names of the triangles of each of the $k$ pentagons surrounding the last circle of $F$ is called similarly $a_2,b_2,c_2$.  So each $\alpha_i$ has the following $3$ neighbors: two $\alpha_i$'s and one $a_i$.  The neighbors of $a_i$ are $\alpha_i$ and two $b_i$'s.  The neighbors of $b_i$ are $a_i,b_i,c_i$.

 Next we divide each hexagon to six triangles of the types $A,B,C$.  The two $A$'s have one edge on a circle $p$ or $p=2,\ldots,\ell-1$.
 Then we have two triangles $B$.  Note that the interior of each edge of $B$ lies in the interior of the hexagon.  Then we have two triangles
 $C$.  Each $C$ have one edge on the circle $p+1$ of $F$.   Assume for simplicity of the argument that $\ell\ge 4$, i.e. we have at least three layer of hexagons. So the neighbors of $c_1$ are $b_1,c_1,A$.   The neighbors of $c_2$ are $b_2,c_2,C$.

 We now analyze the neighbors of the six triangles induced by hexagons.  First, let us consider the first layer of hexagons between the second and the third
 layer of $F$.
 The neighbors of $A$ are $c_1,A,B$.  The neighbors of $B$ are $A,B,C$.  (This fact holds for any hexagon.) The neighbors of $C$ are $B,C$ located inside the
 the same hexagon and a triangle $A$ locates in the next level of hexagons, between  the third and the fourth circles.  (This fact holds for all layers
 of hexagons except the last layer located between the circle $\ell$ and $\ell+1$ of $F$.
 For all other hexagons located between the circle $p$ and $p+1$ of $F$ for $p=3,\ldots,\ell+1$ the neighbors of $A$ are $A,B,C$.
 Finally, in the last layer of hexagons, the neighbors of $C$ are $B,C,c_2$.

 We now discuss the semi-circular structure of $Le(F)$.
 We assume that we have $q$ circles.
 The last circle $q$ of $Le(F)$ is $C_k$, composed of $k$ $\alpha_2$.
 The circle  $q-1$ is $C_{3k}$ of the form $a_2-b_2-b_2-a_2-b_2-b_2-\ldots-b_2-b_2-a_2$.  The circle $q-2$ is $C_{4k}$ of the form
 $c_2-c_2-C-C-c_2-c_2-\ldots-C-C-c_2-c_2$.  The $q-3$ circle is $C_{4k}$ of the form is $B-B-A-A-\ldots-A-A-B-B$, since it uses $4$ triangles in each
 hexagon in the last hexagon level.  The $q-4$ circle of $Le(F)$ is $C_{4k}$ of the form $B-C-C-B-B-C-C-B-\ldots-B$, consisting of triangles obtained from
 the hexagons in the level between the circles $\ell-1$ and $\ell$ of $F$.  The circle $q-5$ of $Le(F)$ is $A-A-C-C-\ldots-A-A$, where the triangles $A$
 are in the previous level and the triangles $C$ are from the next level of hexagons. The $q-6$ circle is $C_{4k}$ of the form $B-C-C-B-B-C-C-B-\ldots-B$.

 Suppose first that $\ell$ is odd.  Then it is straightforward to see that the third circle of $Le(F)$ is $C_{4k}$ of the from $A-A-c_1-c_1-A-A-\ldots-A-A$.
 The second circle of $Le(F)$ is $C_{3k}$ of the form $a_1-b_1-b_1-a_1-\ldots-a_1$.  The first circle of $Le(F)$ is $C_k$ consisting of $\alpha_1$ triangles.
 So the circular structure of $Le(F)$ is $C_k-C_{3k}-C_{4k}-\ldots-C_{4k}-C_{3k}-C_{k}$.  Since the number of circles is $q$ the total number of vertices
 is $4k(q-4)+8k$.  Since $F$ had $(\ell+1)2k$ vertices we deduce that $3(\ell+1)2k=4k(q-2)$.  So $q=\frac{3\ell-1}{2}+4$ circles.  So we have exactly
 $\frac{3\ell-1}{2}$ circles $C_{4k}$.  These arguments apply also to the case $\ell=1,3$.

 Assume that $\ell$ is even.  By the above method the last circle is $C_{k}$.  One before the last circle is a circle $C_{3k}$ of the form $A-B-B-A-\ldots-A$.  Then follows a sequence of circles $C_{4k}$.   Finally the first cycle is $C_{k}$.  Inside the circle one we have $k$ triangles of the form $\alpha_1$ and $a_1$.  So the circles of $Le(F)$
 are $C_{k}-C_{4k}-C_{4k}-\ldots-C_{4k}-C_{3k}-C_k$.  The total number of vertices of $Le(F)$ is $(q-3)4k+3k+2k+k $.  Hence $(q-3)4k+6 k=3(l+1)2k$.
 So $Le(F)$ has $\frac{3l}{2}$ $C_{4k}$ circles.  \qed

 \subsection{Upper matching estimates on semi-circular graphs}\label{sub:upmestcr}
 So we can apply the generalized Hadamard's determinant inequality \eqref{FHineq} to a semi-circular graph $G=(V,E)$ by considering the partitioning set of edges of $V$ to $V_1\cup_{j\in [k-1]} V_{j+1}\setminus V_{j}$ .
 To obtain some better upper inequalities we need several auxiliary results.
 We start with the following lemma which is deduced straightforward.
 \begin{lemma}\label{propscycl}  Let $S(C_n)$ be an orientation of $n$-cycle $C_n=([n],E)$ with edges $1-2-\ldots-n-1$.  Then there exists a diagonal matrix $D$ whose diagonal entries are $\pm 1$ such that $DS(C_n)D$ corresponds to the following orientation $1\to 2\to \ldots\to n$
 of these $n-1$ edges in $C_n$.  Let $T_{n,+},T_{n,-}$ be the skew symmetric matrices corresponding to the orientation $1\to 2\to \ldots\to n$
 and $1 \to n, n\to 1$ of $C_n$ respectively.
 \end{lemma}

 Note that
 \begin{equation}\label{detTn}
 \begin{cases} \det T_{n,+}=\det T_{n,-}=0 & \text{if } n \text{ odd}
 \\
 \det T_{n,+}=4,\quad \det T_{n,-}=0  & \text{if } n \text{ even}
 \end{cases}
 \end{equation}

 Indeed if $n$ is odd then the determinant of any skew symmetric matrix is zero.  Assume that $n$ is even.  Then $\det T_{n,-}=0$ since $T_{n,-}\1=\0$, where $\1=(1,\ldots,1)\trans$.  $T_{n,+}$ represents a pfaffian orientation of $C_n$. Hence $\det T_{n,+}$ is the square
 of the perfect matches in $C_n$, which is $2^2$.

 Let $G=(V,E)$ be a simple graph.  Let $W\subset V, |W|\ge 3$ and assume that the subgraph $G(W)$ is a cycle.  Let $W=\{w_1, \ldots,w_{|W|}\}$
 where the edges of $G(W)$ are $w_1-w_2-\ldots-w_{|W|}-w_1$.  Then $D_c(W)$ is a diagonal matrix of order $|W|$ of the form $\diag(d(w_1),\ldots,
 d(w_{|W|}))$.
 \begin{prop}\label{upmatcircgrh}  Let $G=(V,E)$ be a connected semi-circular graph with $k>1$ as described in Lemma \ref{charcircg}.
 Let $S(G)$ be a pfaffian orientation of $G$.  Assume that $B=-S(G)^2$.  Then
 \begin{equation}\label{upmatcircgrh0}
 \perfmat G\le (\det B[V_1])^{\frac{1}{4}}\prod_{i\in[k-1]} (\det B[V_{i+1}\setminus V_{i}])^{\frac{1}{4}}.
 \end{equation}
 Assume furthermore the following conditions.  First, $G$ does not have $3$ and $4$-cycles.  (In particular, each $m_i:=|V_i\setminus V_{i-1}|> 4$ for $i\in [k]$.) Second, each $u\in V_i\setminus V_{i-1}$ is connected to at most one vertex in $V_{i-1}\setminus V_{i-2}$  for each $i>1$ and at most to one vertex in $V_{i+1}\setminus V_i$ for each $i\in [k]$, where  $V_{k+1}=V_k$.
 Then the number of perfect matchings in $G$ is not more than
 \begin{eqnarray} \label{upmatcircgrh1}
 &&(\det B[V_1])^{\frac{1}{4}}\times\\
 \notag
 &&\prod_{i=2}^{k} (\max(\det (D_c(V_{i}\setminus V_{i-1})-2 I_{m_i}-T_{m_i,+}^2), \det (D_c(V_{i}\setminus V_{i-1})-2 I_{m_i}-T_{m_i,-}^2))
 ^{\frac{1}{4}}.
 \end{eqnarray}
 \end{prop}
 \proof  Clearly \eqref{upmatcircgrh0} follows from \eqref{FHineq}.
 Assume now that $G$ satisfies the additional conditions of the theorem.
 Let $B=-S(G)^2=[b_{uv}]_{u,v\in V}$.  Then the diagonal entries $b_{uu}=d(u)$.  Consider the off-diagonal entries of $B$.
 Note that $b_{uv}=0$ if there is no path of length two between $u$ and $v$ in $G$.  Assume that there is a path of length two between $u$ and $v$
 Since $G$ does not have a $4$-cycle the path between $u$ and $v$ of length two is unique.  Hence $b_{uv}=\pm 1$.
 Take any two points $u,v\in V_{i}\setminus V_{i-1}$.  The assumptions of the theorem yield that there is a path of length two between $u$ and $v$
 is two if and only if the distance between $u$ and $v$ on the cycle $G(V_{i}\setminus V_{i-1})$ is two for $i>1$.
 Hence, $B[V_{i}\setminus V_{i-1}]=D_c(V_{i}\setminus V_{i-1})-2 I_{m_i}-S^2(G(V_{i}\setminus V_{i-1}))$, where the orientation of $G(V_{i}\setminus V_{i-1})$ is induced by the orientation $G$.  Recall that there exists a diagonal matrix $D_i$ with entries $\pm 1$
 such that $D_iS(G(V_{i}\setminus V_{i-1}))D_i\in \{T_{m_{i},+}, T_{m_{i},-}\}$.  Clearly $D_iD_c(V_{i}\setminus V_{i-1})D_i=D_c(V_{i}\setminus V_{i-1})$.  Hence
 \[\det B[V_{i}\setminus V_{i-1}]\le \max(\det (D_c(V_{i}\setminus V_{i-1})-2 I_{m_i}-T_{m_i,+}^2), \det (D_c(V_{i}\setminus V_{i-1})-2 I_{m_i}-T_{m_i,-}^2).\]
 \qed

 We now apply the above result to semi-circular cubic graphs satisfying all the above conditions.  Then $D_c(V_{i}\setminus V_{i-1})-2 I_{m_i}=I_{m_i}$.
 To estimate the right-hand side of \eqref{upmatcircgrh1} we need the following lemma.
 \begin{lemma}\label{detTnest}  For $n\ge 3$
 \begin{equation}\label{detTnest2}
 \begin{cases} \det (I_n-T_{n,+}^2)=\det (I_n-T_{n,-}^2)=((\frac{1+\sqrt{5}}{2})^n + (\frac{1-\sqrt{5}}{2})^n)^2 & \text{if } n \text{ odd}
 \\
 \det (I_n-T_{n,+}^2)=((\frac{1+\sqrt{5}}{2})^n + (\frac{1-\sqrt{5}}{2})^n+2)^2, & \text{if } n \text{ even}\\
 \det (I_n-T_{n,-}^2)=((\frac{1+\sqrt{5}}{2})^n + (\frac{1-\sqrt{5}}{2})^n-2)^2.
 & \text{if } n \text{ even}
 \end{cases}
 \end{equation}
 Furthermore, for $n\ge 3$ the subsequence $(\det(I_n-T_{n,+}^2))^{\frac{1}{n}}$ is an increasing subsequence for $n=3,5,7,\ldots$
 and a decreasing subsequence for $n=4,6,8,\ldots$.  Both subsequences converge to $(\frac{1+\sqrt{5}}{2})^2\approx 2.6180$.  Hence
 \begin{equation}\label{detTnest3}
 (\det(I_n-T_{n,+}^2))^{\frac{1}{n}}\le (20)^{\frac{1}{3}}\approx 2.7144 \text { for } n\ge 5.
 \end{equation}
 Equality holds if and only if $n=6$.

 \end{lemma}
 \proof
 Let $S$ be a skew symmetric matrix and $a>0$.  Then
 \begin{eqnarray*}
 &&\det (a^2 I-S^2)=\det((aI-S)(aI+S))=\det (aI+S)\det(aI-S)\\
 &&=\det(aI+S)\det((aI-S)\trans)=\det(aI+S)^2.
 \end{eqnarray*}
 Hence
 \begin{equation}\label{squareid}
 \det (I_n-T_{n,+}^2)=(\det(I_n+T_{n,+}))^2, \quad \det (I_n+T_{n,-}^2)=(\det(I_n-T_{n,-}))^2.
 \end{equation}
 Next we claim that
 \begin{equation}\label{dettT-}
 \det(I_n+tT_{n,-})=\sum_{j=0}^{\lceil\frac{n}{2}\rceil-1} \phi(j,C_n)t^{2j}.
 \end{equation}

 First observe that since $\det T_{n,+}=0$ the polynomial $\det(I_n+tT_{n,-})$ is of degree $n-1$ at most.
 Let $k\in [n-1]$.  Then the coefficient of $t^k$ is the sum of all $\det T_{n,-}[J]$ where $J$
 is a a subset of cardinality of $k$ of $[n]$.  As $T_{n,-}[J]$ represents an orientation of a subgraph $C_n[J]$, which is a forest.
 As any orientation of the forest is pfaffian we deduce  that $\det C_n[J]$ is a square of the number of perfect matches in the subgraph $C_n[J]$.
 Observe next a forest has either $0$ or $1$ perfect match.  Hence $\det C_n[J]$ is the number of perfect matches in the subgraph $C_n[J]$.
 Thus $\sum_{|J|=k}\det T_{n,-}[J]$ is zero if $k$ is odd, and $\phi(\frac{k}{2},C_n)$ if $k$ is even.  Since $\det I_n=1$ we deduce
 \eqref{dettT-}.

 Next we claim that
 \begin{equation}\label{dettT+}
 \begin{cases}
 \det (I_n+tT_{n,+})=\det (I_n+tT_{n,-}) & \text{if } n \text{ odd}
 \\
 \det (I_n+tT_{n,+})=\det (I_n+tT_{n,-})+4t^n & \text{if } n \text{ even}
 \end{cases}
 \end{equation}

 Indeed, the above arguments show that $\det (I_n+tT_{n,+})-\det (I_n+tT_{n,-})=t^n\det T_{n,+}$.  Hence \eqref{dettT+} holds.
 Let $\phi(t,C_n):=\sum_{j=0}^{\lfloor\frac{n}{2}\rfloor}\phi(j,C_n) t^{2j}$ be the matching polynomial of $C_n$.
 Combine \eqref{dettT-} and  \eqref{dettT+} to deduce
 \begin{equation}\label{detT+matpol}
 \begin{cases}
 \det(I_n+tT_{n,+})=\phi(t,C_n) & \text{if } n \text{ odd}\\
 \det(I_n+tT_{n,+})=\phi(t,C_n)+2t^n & \text{if } n \text{ even}
 \end{cases}
 \end{equation}

 Let $P_n$ be a path of length $n$.  Denote by $\phi(t,P_n)$ the matching polynomial of a path of length $n$.
 Then the following identities are straightforward, see \cite[(2.5),(2.6)]{FKM}:
 \begin{eqnarray}\label{relPcons}
 &&\phi(t,P_n)=\phi(t,P_{n-1})+t\phi(t,P_{n-2}),\\
 &&\phi(t,C_n)=\phi(t,P_n)+t\phi(t,P_{n-2}).\label{relPC}
 \end{eqnarray}
 It easily follows that $\phi(1,P_n)$ is the standard Fibonacci sequence for $1,1,2,3,5,\ldots$ for $n=0,1,2,3,4,\ldots$.
 Use this fact and \eqref{relPC} to obtain
 \begin{equation}\label{totmatCn}
 \phi(1,C_n)=(\frac{1+\sqrt{5}}{2})^n + (\frac{1-\sqrt{5}}{2})^n \text{ for } n\ge 3.
 \end{equation}

 Combine \eqref{totmatCn}, \eqref{detT+matpol}, \eqref{dettT+} and \eqref{squareid} to deduce \eqref{detTnest2}.
 Denote $\epsilon:=\frac{\sqrt{5}-1}{\sqrt{5}+1}, \delta:=\frac{2}{\sqrt{5}+1}$.
 Let
 \[a_n:=(\det(I_n-T_{n,+}^2))^{\frac{1}{n}}, b_n:=(\frac{\sqrt{5}+1}{2})^2( 1-\epsilon^n)^{\frac{2}{n}}, c_n:=
 (\frac{\sqrt{5}+1}{2})^2( 1+\epsilon^n+2\delta^n)^{\frac{2}{n}},\;n\in\N.\]
 Then this subsequence is equal to the subsequence  for $n=3,5,7,\ldots$.
 Then $a_n=b_n$ for odd $n$ and $a_n=c_n$ for an even $n$.  As $\epsilon,\delta \in (0,1)$ we deduce
 \begin{eqnarray*}
 &&(1-\epsilon^{n})^{\frac{1}{n}}<(1-\epsilon^{n+1})^{\frac{1}{n}}<(1-\epsilon^{n+1})^{\frac{1}{n+1}},\\
 &&( 1+\epsilon^n+2\delta^n)^{\frac{1}{n}} > ( 1+\epsilon^{n+1}+2\delta^{n+1})^{\frac{1}{n}} > ( 1+\epsilon^{n+1}+2\delta^{n+1})^{\frac{1}{n+1}}.
 \end{eqnarray*}

 Hence $a_n, n=1,3,\ldots$ is an increasing sequence and $a_n, n=2,4,\ldots,$ is a decreasing sequence which both converge to
 $(\frac{\sqrt{5}+1}{2})^2$.  Therefore inequality \eqref{detTnest3} holds.  \qed
 \begin{theo}\label{upbdssemcrc}  Let $G=(V,E), n=|E|$ be a semi-circular cubic graph where $n_1=|V_1|$ and $n_0=|V_0|$ is the number of vertices of $G$ in the the closed disk bounded by the first circle $O_1$ and its interior respectively.  If $G$ does not have $3$ and $4$-cycles then
 \begin{eqnarray}\label{upbdssemcrc1}
 &\perfmat G\le 20^{\frac{n-n_1}{12}}3^{\frac{n_1}{4}} & \text{ if } n_0>0,\\
 \label{upbdssemcrc2}
 &\perfmat G\le 20^{\frac{n}{12}} & \text{ if G is circular, i.e. } n_0=0.
 \end{eqnarray}
 \end{theo}
 \proof
 Note that if $G$ is a semi-circular cubic graph then each vertex on a circle $O_i$ for $i>1$ has exactly one edge which connects it either to $O_{i-1}$ or to $O_{i+1}$.  Hence \eqref{upmatcircgrh1} holds.  Clearly $\det B[V_1]\le 3^{\frac{n_1}{3}}$. Use Lemma \ref{detTnest} to deduce \eqref{upbdssemcrc1}.
 Assume now that $n_0=0$, i.e. $G$ is circular.  Then the arguments of the proof of Proposition \ref{upmatcircgrh} yields that
 $\det B[V_1]\le \max(\det (I_{m_1}+T_{m_1,+}^2),\det (I_{m_1}+T_{m_1,-}^2))$, where $|V_1|=n_1=m_1$.  Lemma \ref{detTnest}
 yields that $\det B[V_1]\le 20^{\frac{m_1}{3}}$.  Hence \eqref{upbdssemcrc2} holds.  \qed

 Since $20^{\frac{1}{3}}\approx 2.7144 < \sqrt{8}\approx 2.8284$ we deduce that \eqref{upbdssemcrc2} is better than the best
 estimate $8^{\frac{n}{8}}$ for fullerene graphs with hamiltonian cycle given in \eqref{estM'a}.

 In applying Theorem \ref{upbdssemcrc} to fullerene graphs one needs to have a good estimate on $n_1$.  The examples that we considered
 suggest that there exists a universal constant $N_1$ such that any fullerene graph is a semi-circular graph with $|V_1|\le N_1$.

 The next obvious problem is how good is the upper bound \eqref{upbdssemcrc2} for circular fullerene graph.
 Recall that a fullerene with nontrivial cyclic 5-edge cut one has at least $5^{\frac{n-20}{10}}$ perfect matches \cite{QZ05} and
  \cite[8-4]{LoeblL8}.

\end{document}